\documentclass[12pt]{article}
\usepackage{graphicx}
\usepackage{amsmath}
\usepackage{amssymb}
\usepackage{amscd}
\usepackage{amsthm}
\usepackage{epsfig}
\usepackage{latexsym}
\usepackage{enumerate}

\newtheorem{theorem}{Theorem}[section]
\newtheorem{lemma}[theorem]{Lemma}

\newtheorem{definition}[theorem]{Definition}

\newtheorem{remark}[theorem]{Remark}

\numberwithin{equation}{section} \numberwithin{theorem}{section}

\title{Spectral Scales and Linear Pencils.}
\author{Christopher M. Pavone}
\date{November 07, 2005}

\newcommand{\M}{\mathcal{M}}

\begin{document}

\bibliographystyle{amsplain}
\maketitle

\begin{abstract}
\noindent \emph{Developed in 1999 by Akemann, Anderson, and
Weaver, the spectral scale of an $n\times n$ matrix $A$,
 is a convex,
compact subset of $\mathbb{R}^3$ that reveals important spectral
information about $A$ \cite{AAW}.  In this paper we present new
information found in the spectral scale of a matrix. Given a matrix
$A=A_1 + iA_2$ with $A_1$ and $A_2$ self-adjoint and $A_2\neq 0,$ we
show that faces in the boundary of the spectral scale of $A$ that
are parallel to the $x$-axis describe elements of
$\sigma(A_1,A_2)\bigcap\mathbb{R},$ the real elements of the
spectrum of the linear pencil $P(\lambda)=A_1 + \lambda A_2.$}

\end{abstract}

\section{Preliminary Definitions}

\noindent\textbf{Notation.}  Let $\M = M_n(\mathbb{C}),$ the set
of $n\times n$ matrices with complex entries, and let
$<\cdot,\cdot>$ denote the standard inner-product on
$\mathbb{C}^n.$ For each $C \in \M,$ we let $\tau(C) =
\frac{1}{n}tr(C),$ where $tr$ is the standard canonical trace on
$\M.$ Given $A\in \M,$ $A^*$ will denote the adjoint of $A$,
$\sigma(A)$ will denote the set of eigenvalues (spectrum) of $A$,
and $W(A)$ will denote the numerical range of $A$. We let $||A||$
 denote the operator (spectral) norm of $A$.

\begin{definition}
\emph{Given $A \in \M,$ write $A = A_1 +iA_2$ with $A_1$ and $A_2$
self-adjoint. The \textbf{spectral scale} of $A$ is defined as
\[ B(A) = \left\{(\tau(C),\tau(A_1C),\tau(A_2C)) \ : \ C \in
\M_1^+\right\} \subset \mathbb{R}^3,\] where $\M_1^+$ denotes the
positive semidefinite part of the unit ball of $\M.$}
\end{definition}

\begin{remark}
\emph{$B(A)$ is a convex, compact, subset of $\mathbb{R}^3.$  In the
case $A=A^*$, $B(A)$ determines $A$ up to unitary equivalence and
gives a complete description of $\sigma(A), tr(A), ||A||, \mbox{ and
} W(A).$ When $A\neq A^*$ $B(A)$ completely describes $W(A)$, and in
some cases (e.g., when $A$ is normal), also yields a complete
description of $\sigma(A)$. In its most general setting, the
spectral scale was originally developed for operators that lie in a
finite von Neumann algebra equipped with a finite, faithful, normal
trace \cite{AAW}. We note that $\M = M_n(\mathbb{C})$ is just one
particular setting where the spectral scale can be computed. For the
complete series on the spectral scale see \cite{AAW,AA1,AA2,AA3}.}
\end{remark}

\begin{definition} \emph{If $\mathfrak{F}$ is a nonempty, convex subset of a convex set
$\mathfrak{C}$, then $\mathfrak{F}$ is a \textbf{face} of
$\mathfrak{C}$ if $\alpha x + (1-\alpha)y \in \mathfrak{F}$ implies
$x,y\in \mathfrak{F}$ for all $x,y \in \mathfrak{C}$ and $0 < \alpha
< 1.$}
\end{definition}

\begin{definition}
\emph{Let $A_1, A_2 \in \M$. The matrix valued function
$P:\mathbb{C} \longrightarrow \M$ defined by \[P(\lambda) = A_1 +
\lambda A_2, \ \  \lambda \in \mathbb{C}\] is called a
\textbf{linear (matrix) pencil}.  If $A_1=A_1^*$ and $A_2=A_2^*,$
then $P(\lambda)$ is said to be \textbf{self-adjoint.} }
\end{definition}

\begin{definition}
\emph{The \textbf{spectrum of the pencil} $P(\lambda) = A_1 +
\lambda A_2$, denoted  by $\sigma(A_1,A_2)$, is the subset of
$\mathbb{C}_{\infty} = \mathbb{C}\cup\{\infty\}$ determined by the
following properties. The point $\infty \in \sigma(A_1,A_2)$ if
and only if $A_2$ is not invertible, and
 \[\sigma(A_1,A_2)\cap\mathbb{C} = \{ \lambda \in \mathbb{C} \ : \ A_1 + \lambda
A_2 \mbox{ is not invertible }\}.\] It can be shown that
$\sigma(A_1,A_2)$ is compact (with respect to the usual topology
on $\mathbb{C}_{\infty}$) (\cite{Go}, Section IV.1).}
\end{definition}

\begin{remark}
\emph{Spectral problems for matrix pencils (and more generally
operator polynomials) arise in different areas of mathematical
physics like differential equations, boundary value problems,
controllable systems, the theory of oscillations and waves,
elasticity theory, and hydromechanics \cite{M}.  For more on linear
pencils see \cite{Go,GLR, M}}.

\end{remark}

\noindent\textbf{Notation.}  For the remainder of this we paper fix
$A = A_1 + iA_2$ in $\M$ with $A_j=A_j^*$ ($j=1,2$) and $A_2\neq 0.$
Also, for each $\textbf{t} = (t_1,t_2) \in \mathbb{R}^2$ we form the
new self-adjoint operator $A_{\bf{t}} = t_1A_1 +t_2A_2$. Note that
if we fix $t_1\neq 0$ and let $t_2$ vary, then
\[A_{\bf{t}} = t_1A_1 +t_2A_2 = t_1(A_1 + \frac{t_2}{t_1}A_2)\] is
a multiple of the (real) pencil
$P_{|_{\mathbb{R}}}:\mathbb{R}\longrightarrow \M$ defined by $P(r)
= A_1 + rA_2.$

Given $\textbf{t} = (t_1,t_2) \in \mathbb{R}^2$ with
$t_1^2+t_2^2=1,$  define
\begin{equation}
\theta_{\textbf{t}}=\tan^{-1}\left(\frac{t_2}{t_1}\right)\ \
-\frac{\pi}{2}<\theta_{\textbf{t}}<\frac{\pi}{2},
\end{equation} let

\begin{equation}Q_{\textbf{t}} = \left(%
\begin{array}{ccc}
  1 & 0 & 0 \\
  0 & t_1^2 & t_1t_2 \\
  0 & t_1t_2 & t_2^2 \\
\end{array}%
\right)\end{equation}  be the orthogonal projection of
$\mathbb{R}^3$ onto $span\{(1,0,0),
(0,t_1,t_2)\},$ and let \begin{equation}R_{\textbf{t}} = \left(%
\begin{array}{ccc}
  1 & 0 & 0 \\
  0 & t_1 & -t_2 \\
  0 & t_2 & t_1 \\
\end{array}%
\right) = \left(%
\begin{array}{ccc}
  1 & 0 & 0 \\
  0 & \cos(\theta_{\bf{t}}) & -\sin(\theta_{\bf{t}}) \\
  0 & \sin(\theta_{\bf{t}}) & \cos(\theta_{\bf{t}}) \\
\end{array}%
\right)\end{equation} denote the rotation of $\mathbb{R}^3$ about
the $x$-axis through an angle of $\theta_{\textbf{t}}$ (see Figure
1).

\section{A New Result for $B(A)$}

We start by describing the relationship between $B(A)$ and
$B(A_{\bf{t}})$.

\begin{figure}
\ \ \ \ \ \ \ \ \ \ \ \ \ \ \ \ \ \ \ \ \ \ \
\includegraphics[scale=1.2]{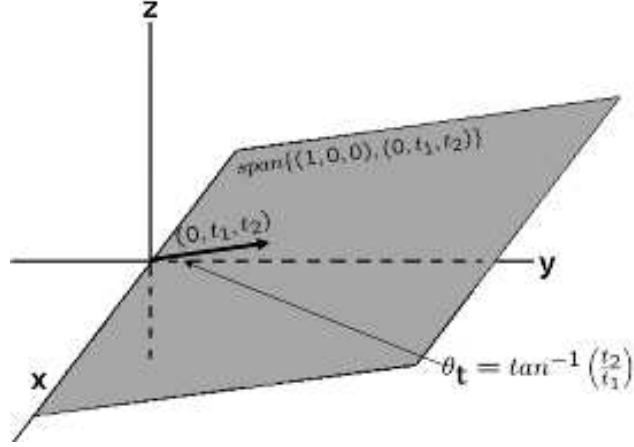}
\caption{$\theta_{\textbf{t}}=\tan^{-1}\left(\frac{t_2}{t_1}\right)$}
\end{figure}

\begin{theorem}
Let $\textbf{t} = (t_1,t_2) \in \mathbb{R}^2$ such that
$t_1^2+t_2^2=1$.   Then
$Q_{\textbf{t}}(B(A))=R_{\textbf{t}}(B(A_{\bf{t}})).$
\end{theorem}

\begin{proof}
Choose $\textbf{t} \in \mathbb{R}^2$  as above, and define
$\theta_{\textbf{t}},$ $Q_{\textbf{t}}$ and $R_{\textbf{t}}$ as in
(1.1), (1.2), and (1.3) respectively.   If
$v=(\tau(C),\tau(A_1C),\tau(A_2C))$ is any point in $B(A),$ then
\begin{align*} Q_{\bf{t}}(v) &= \left(\tau(C), t_1^2\tau(A_1C) +
t_1t_2\tau(A_2C), t_1t_2\tau(A_1C) + t_2^2\tau(A_2C)\right)\\ &=
\left(\tau(C), t_1(\tau(A_{\textbf{t}}C)),
t_2(\tau(A_{\textbf{t}}C)\right)\\
&=\left(\tau(C),\cos(\theta_{\textbf{t}})\tau(A_{\textbf{t}}C),\sin(\theta_{\textbf{t}})\tau(A_{\textbf{t}}C)\right)\\
&=\left(%
\begin{array}{ccc}
  1 & 0 & 0 \\
  0 & \cos(\theta_{\bf{t}}) & -\sin(\theta_{\bf{t}}) \\
  0 & \sin(\theta_{\bf{t}}) & \cos(\theta_{\bf{t}}) \\
\end{array}%
\right)\left(%
\begin{array}{c}
  \tau(C) \\
  \tau(A_{\textbf{t}}C) \\
  0 \\
\end{array}%
\right)\\
&=R_{\textbf{t}}((\tau(C),\tau(A_{\textbf{t}}C),0)).\end{align*}
From this we see that $Q_{\bf{t}}(B(A)) =
R_{\textbf{t}}(B(A_{\textbf{t}})).$ Thus, the orthogonal projection
of the spectral scale of $A$ onto $span\{(1,0,0), (0,t_1,t_2)\}$, is
a rotation of the spectral scale of $A_{\textbf{t}}$ about the
$x$-axis, through an angle of $\theta_{\bf{t}}$.
\end{proof}

\begin{remark}
\emph{This result is different from that in (\cite{AAW}, Lemma
2.2). There it is shown that if $\pi_{\textbf{t}} = \left(%
\begin{array}{ccc}
  1 & 0 & 0 \\
  0 & t_1 & t_2 \\
  0 & 0 & 0 \\
\end{array}%
\right),$  then $\pi_{\textbf{t}}(B(A)) = B(A_{\textbf{t}}).$
Although this is an easy way to extract $B(A_{\bf{t}})$ from
$B(A)$, it is not as nice geometrically, as Figure 2 demonstrates.
Furthermore, (\cite{AAW}, Lemma 2.2) puts no restrictions on
$\textbf{t}=(t_1,t_2).$  If the restriction $t_1^2+t_2^2=1$ is put
on \textbf{t}, then we see that $\pi_{\textbf{t}} =
R_{\textbf{t}}^*Q_{\textbf{t}}$. }
\end{remark}

\begin{figure}
\ \ \ \ \ \ \ \ \ \ \ \ \ \ \ \ \ \
\includegraphics[scale=0.5]{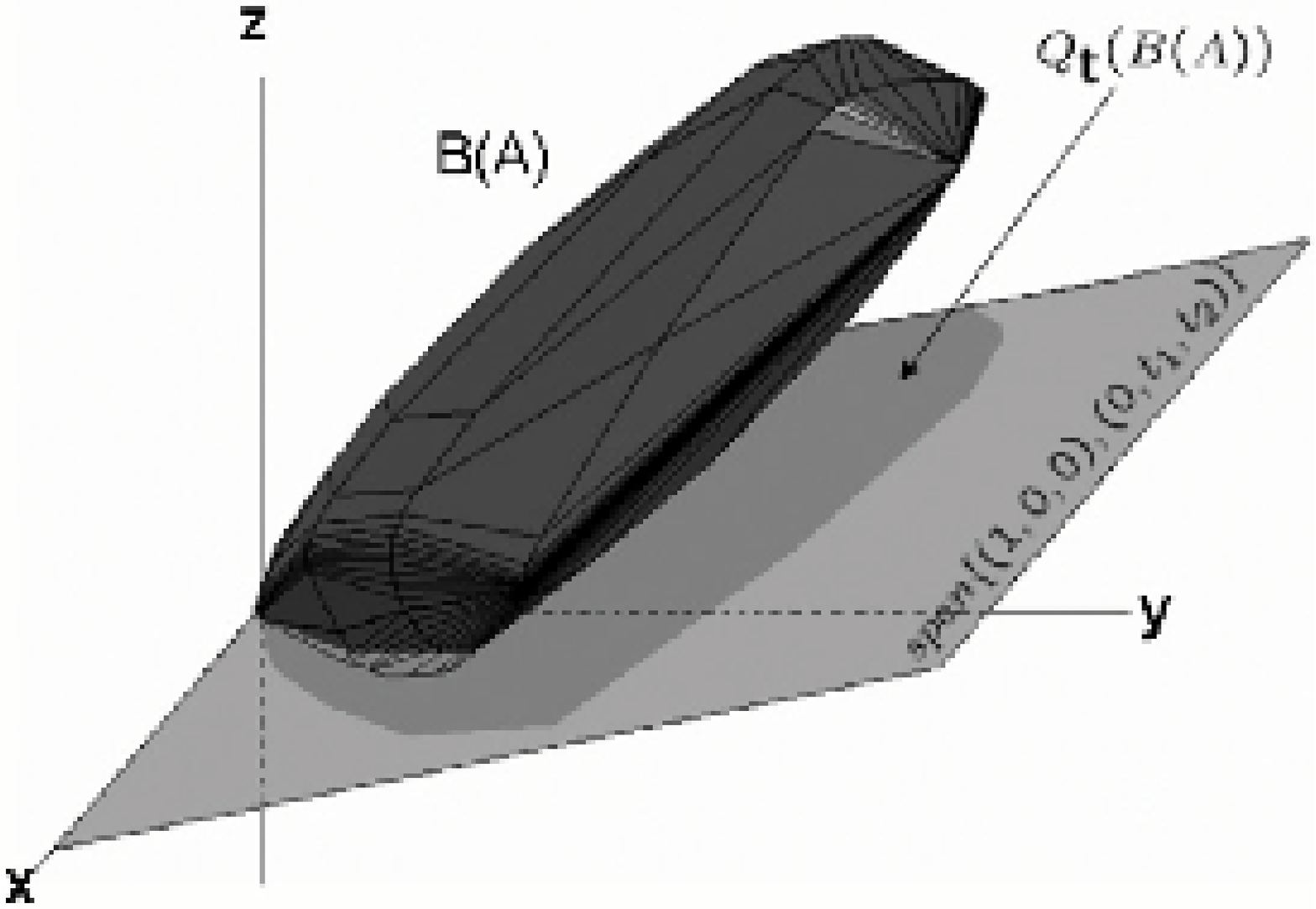}
\caption{$Q_{\bf{t}}(B(A)) = R_{\textbf{t}}(B(A_{\textbf{t}}))$}

\end{figure}

\begin{lemma}
For $\frac{\pi}{2} \leq \theta_{\textbf{t}} \leq \frac{\pi}{2},$
$\tan(\theta_{\textbf{t}}) \in \sigma(A_1,A_2)$ iff $0 \in
\sigma(A_{\bf{t}}).$
\end{lemma}

\begin{proof}  For $-\frac{\pi}{2} < \theta_{\textbf{t}} <
\frac{\pi}{2},$  $A_1 + \tan(\theta_{\textbf{t}})A_2$ is not
invertible iff $\cos(\theta_{\textbf{t}})A_1 +
\sin(\theta_{\textbf{t}})A_2$ is not invertible iff $0 \in
\sigma(\cos(\theta_{\textbf{t}})A_1 +
\sin(\theta_{\textbf{t}})A_2) = \sigma(t_1A_1 + t_2A_2) =
\sigma(A_{\textbf{t}}).$ For $\theta_{\textbf{t}} =
\pm\frac{\pi}{2},$ we have $\textbf{t}=(0,\pm 1)$ and
$\tan(\theta_{\textbf{t}}) = \infty \in \sigma(A_1,A_2)$ iff $A_2$
is not invertible iff $0\in\sigma(A_2) = \sigma(A_{\textbf{t}}),$
and our claim holds.
\end{proof}

\begin{remark}  \emph{Our final theorem (Theorem 2.5) tells us how to read
$\sigma(A_1,A_2)\bigcap \mathbb{R}$ from $B(A).$ Note that if $A$
is normal and $\lambda = x+iy \in \mathbb{C}$ is not real, then if
$(A_1 + \lambda A_2)\eta = 0$ for some nonzero
$\eta\in\mathbb{C}^n$, then we must have \begin{align*}0 &= <(A_1
+ \lambda
A_2)\eta,A_2\eta>\\
&=<A_1\eta,A_2\eta> + \lambda<A_2\eta,A_2\eta>\\
&=<A_1\eta,A_2\eta> + x<A_2\eta,A_2\eta> + iy<A_2\eta,A_2\eta>\\
&=(<A_2A_1\eta,\eta> + x<A_2\eta,A_2\eta>) + i(y<A_2\eta,A_2\eta>).
\end{align*} Since $<A_2A_1\eta,\eta> \in \mathbb{R},$ we get that $A_2\eta=0$ (i.e., $A_2$ is not invertible).  Thus,
when $A$ is normal and $A_2$ is invertible we have $\sigma(A_1,A_2)
\subset \mathbb{R}\cup\{\infty\}.$ In this case, Theorem 2.5 gives
us a complete description of the spectrum of the self-adjoint pencil
$P(\lambda)=A_1+\lambda A_2.$}
\end{remark}

\begin{theorem}  Let $\textbf{t}=(t_1,t_2) \in \mathbb{R}^2$ such that $t_1^2+t_2^2=1.$  The following statements are equivalent:

\begin{enumerate}

\item $0 \in \sigma(A_{\textbf{t}}).$

\item $\tan(\theta_{\textbf{t}})\in\sigma(A_1,A_2).$

\item There exists a 1-dimensional face
$\mathfrak{F}_{\textbf{t}}$ in the boundary of $B(A_{\textbf{t}})$
that is parallel to the $x$-axis.

\item There exists a 1-dimensional face
$\mathfrak{F}_{\textbf{t}}'$ in the boundary of
$Q_{\textbf{t}}(B(A))$ that is parallel to the $x$-axis.

\item There exists a 1-dimensional (or 2-dimensional) face
$\mathfrak{F}_{\textbf{t}}''$ in the boundary of $B(A)$ that is
parallel to the $x$-axis.

\end{enumerate}
\end{theorem}

\begin{proof} By Lemma 2.3 we have $\tan(\theta_{\textbf{t}})\in\sigma(A_1,A_2)$ iff
$0\in \sigma(A_{\textbf{t}}),$ so (1)$\Leftrightarrow$(2). By
(\cite{AAW}, Theorem 1.5), this is equivalent to the existence of a
1-dimensional face $\mathfrak{F}_{\textbf{t}}$ contained in the
boundary of $B(A_{\textbf{t}})$ which is parallel to the $x$-axis,
thus (1)$\Leftrightarrow$(3). Clearly, if a line segment is parallel
to the $x$-axis, then a rotation about the $x$-axis does not change
this fact.  Therefore, by Theorem 2.1,
$\mathfrak{F}_{\textbf{t}}'=R_{\textbf{t}}(\mathfrak{F}_{\textbf{t}})$
is a 1-dimensional face contained in the boundary of
$Q_{\textbf{t}}(B(A_{\textbf{t}}))$ which is parallel to the
$x$-axis, so (3)$\Leftrightarrow$(4). But then, since
$Q_{\textbf{t}}$ is linear we have that
$\mathfrak{F}_{\textbf{t}}''=Q_{\textbf{t}}^{-1}(\mathfrak{F}_{\textbf{t}}')$
is a face of $B(A),$ and it is also parallel to the $x$-axis (see
Figure 3).  Therefore (4)$\Leftrightarrow$(5), and our claim holds.
\end{proof}

\begin{figure}
\ \ \ \ \ \ \ \ \ \ \ \ \ \
\includegraphics[scale=1.0]{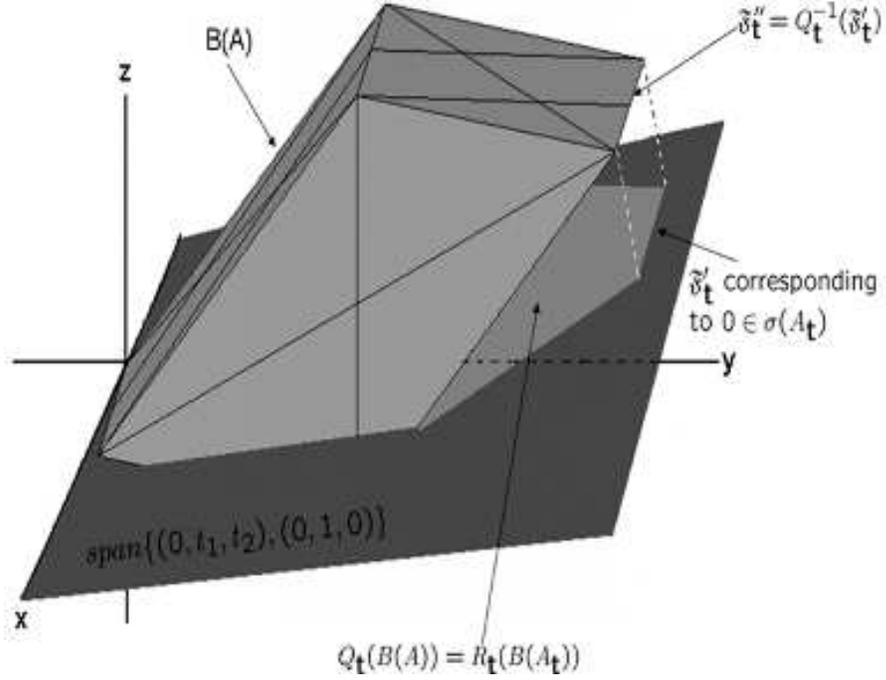} \\
  \caption{Theorem 2.5}
\end{figure}

\newpage

\bibliography{pencils}

\providecommand{\bysame}{\leavevmode\hbox to3em{\hrulefill}\thinspace}
\providecommand{\MR}{\relax\ifhmode\unskip\space\fi MR }
\providecommand{\MRhref}[2]{%
  \href{http://www.ams.org/mathscinet-getitem?mr=#1}{#2}
}
\providecommand{\href}[2]{#2}
\begin{thebibliography}{1}

\bibitem{AA3}
Charles~A. Akemann and Joel Anderson, \emph{A geometric spectral theory for
  {$n$}-tuples of self-adjoint operators in finite von {N}eumann algebras.
  {II}}, Pacific J. Math. \textbf{205} (2002), no.~2, 257--285.

\bibitem{AA1}
\bysame, \emph{The spectral scale and the {$k$}-numerical range}, Glasg. Math.
  J. \textbf{45} (2003), no.~2, 225--238.

\bibitem{AA2}
\bysame, \emph{The spectral scale and the numerical range}, Internat. J. Math.
  \textbf{14} (2003), no.~2, 171--189.

\bibitem{AAW}
Charles~A. Akemann, Joel Anderson, and Nik Weaver, \emph{A geometric spectral
  theory for {$n$}-tuples of self-adjoint operators in finite von {N}eumann
  algebras}, J. Funct. Anal. \textbf{165} (1999), no.~2, 258--292.

\bibitem{Go}
Goldberg S. Kaashoek~M.A. Gohberg, I., \emph{Classes of linear operators
  vol.1}, Birkhauser Verlag, Germany, 1990.

\bibitem{GLR}
Lancaster P. Rodman~L. Gohberg, I., \emph{Matrix polynomials}, Academic Press,
  New York, 1982.

\bibitem{M}
A.S. Markus, \emph{Introduction to the spectral theory of polynomial operator
  pencils}, Translations of Mathematical Monographs, vol.~71, American
  Mathematical Society, Providence, RI, 1988.

\end{thebibliography}

\end{document}